\documentclass{article}
\usepackage{graphicx}
\usepackage{amsmath}
 \usepackage{amsfonts}
 \usepackage{color}

\begin{document}

\title{Non-Power Positional Number Representation Systems, Bijective Numeration, and the Mesoamerican Discovery of Zero} 

\author{
Berenice Rojo-Garibaldi${}^{a}$, Costanza Rangoni${}^{b}$, \\
Diego L. Gonz\'alez${}^{b,c}$,  and Julyan H. E. Cartwright${}^{d,e}$ \\
${}^{a}$ Posgrado en Ciencias del Mar y Limnolog\'{\i}a, \\ Universidad Nacional Aut\'onoma de M\'exico, \\
Av.\ Universidad 3000, Col.\ Copilco, \\ Del.\ Coyoac\'an, Cd.Mx.\ 04510, M\'exico\\
${}^{b}$
Istituto per la Microelettronica e i Microsistemi, \\ Area della Ricerca CNR di Bologna, 40129 Bologna,  Italy\\
${}^{c}$
Dipartimento di Scienze Statistiche ``Paolo Fortunati'', \\ Universit\`a di Bologna, 40126 Bologna, Italy\\
${}^{d}$
Instituto Andaluz de Ciencias de la Tierra, \\ CSIC--Universidad de Granada, 18100 Armilla, Granada, Spain\\
${}^{e}$
Instituto Carlos I de F\'{\i}sica Te\'orica y Computacional, \\ Universidad de Granada, 18071 Granada, Spain
}

\date{}

\maketitle


\begin{abstract}
Pre-Columbian Mesoamerica was a fertile crescent  for the development of number systems. A form of vigesimal system seems to have been present from the first Olmec civilization onwards, to which succeeding peoples made contributions.
We discuss the Maya use of the representational redundancy present in their Long Count calendar, a non-power positional number representation system with multipliers 1, 20, 18$\times$ 20, $\ldots$, 18$\times$ 20$^n$. We demonstrate that the Mesoamericans did not need to invent positional notation and discover zero at the same time because they were not afraid of using  a number system 
in which the same number can be written in different ways.
A Long Count number system with digits from 0 to 20 is seen later to pass to one using digits 0 to 19, which leads us to propose that even earlier there may have been an initial zeroless bijective numeration system whose digits ran from 1 to 20.
Mesoamerica was able to make this conceptual leap to the concept of a cardinal zero to perform arithmetic owing to a familiarity with multiple and redundant number representation systems.
\end{abstract}

\section{Introduction}

\begin{figure}
  \centering\includegraphics[width=0.49 \columnwidth]{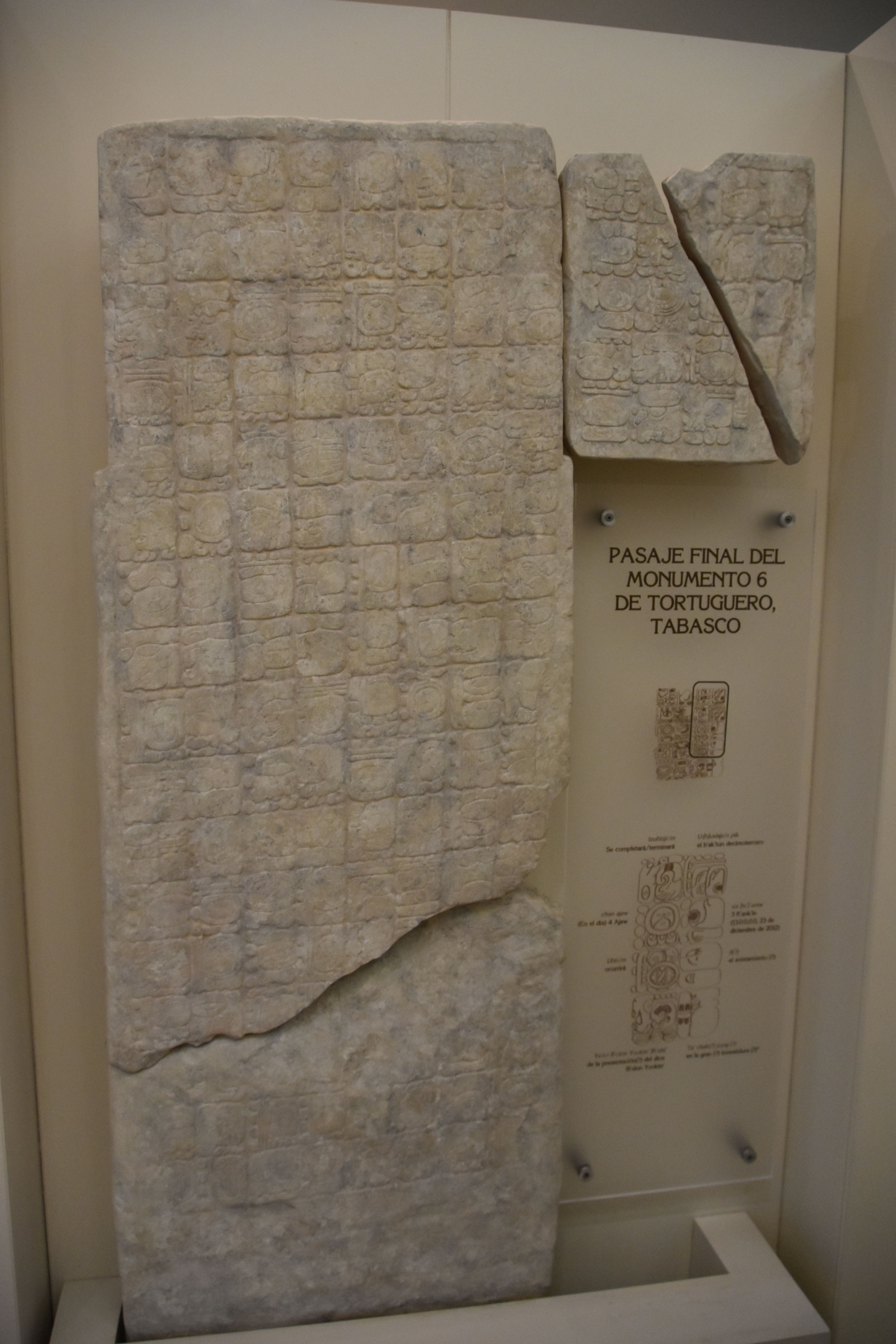} \includegraphics[width=0.49 \columnwidth]{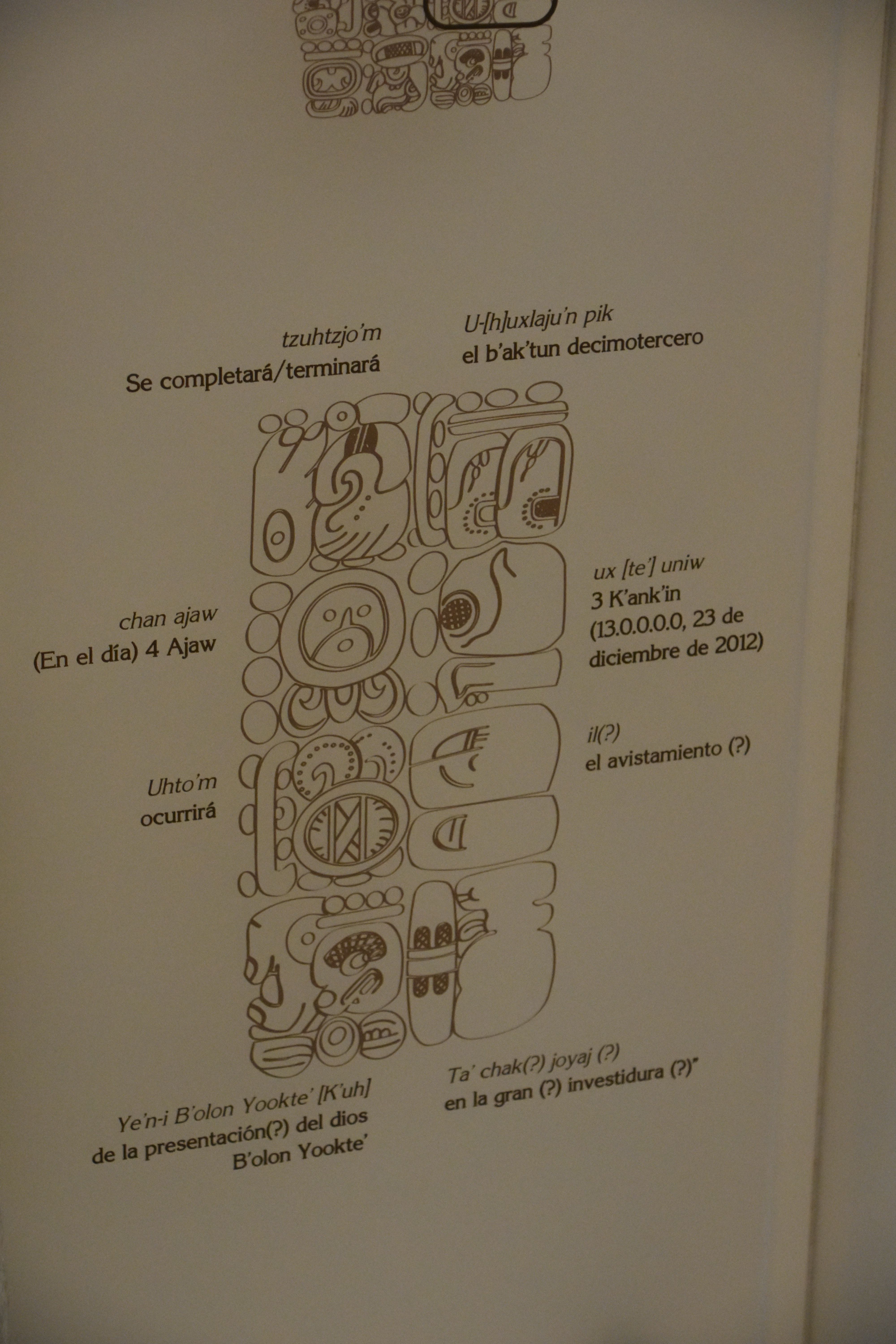}
\caption{The Maya Long Count calendar. The stela on the left contains the dates shown enlarged on the right, including the date of the last completion of a full cycle of the Maya calendar, which occurred on 13.0.0.0.0, 23rd December 2012, many centuries in the future when this stela was inscribed. (Museo Maya de Canc\'un, Instituto Nacional de Antropolog\'{\i}a e Historia, M\'exico.)
\label{fig:Maya}}
     \end{figure}

Alongside the decimal positional number system, Fibonacci popularized a new number in Europe: zero. Partially owing to this historical link, it has almost been, as it were, a truth universally acknowledged that a civilization in possession of a good number system must be in want of a zero. But this is not necessarily so. It is perfectly possible to have a positional number system without a zero (see the Methods section for an introduction to number representation systems). This is called bijective numeration, and we argue that Mesoamerica may well have invented the positional number system first as a bijective system without a zero. Only some time later do we see zeros beginning to appear in the Maya Long Count, depicted in Fig.~\ref{fig:Maya}. Because the Maya were used to a redundant number system they were not afraid of writing the same number in various ways, and they found that the zero they had discovered and initially used in a non-positional system could be introduced into their positional system with minimal problems. Thus they were able to make the conceptual leap to a cardinal zero --- a zero used in arithmetic --- in stages aided by their familiarity with multiple number representation systems.

The organization of this work is as follows: Section 2 introduces the Maya calendar; Section 3 introduces the necessary mathematics; and Section 4 contains the results. The central result is encapsulated in Figure 4c. Section 5 gives our conclusions and Section 6   provides details on the mathematics involved.

\section{The Maya Long Count Calendar}

\begin{quote}
``No people in history has shown such interest in time as the Maya. Records of its passage were inscribed on practically every stela, on lintels of wood and stone, on stairways, cornices, friezes and panels'' J. E. S. Thompson \cite{thompson}.
\end{quote}

\noindent
The Maya \cite{sharer2005,mckillop2006,coe2011} understood well what we now call deep time. 
The Long Count is a positional notation system that, as its name indicates, enables complicated arithmetical calculations over arbitrarily long periods of time. By the classical period as it reached its apogee under the Maya, scribes were writing of time periods of millions of years into the past and thousands of years into the future \cite{blume}.

The Maya developed a very sophisticated astronomical  culture \cite{bricker} in their civilization centred around the Yucat\'an peninsula in what is today Mexico, Guatemala and Belize, whose classic period of greatest splendour ran from around 3rd century  to 10th century CE before falling into decline \cite{turner2012}. 
A numerical calendar is a revolutionary idea: to enumerate the passage of time, rather than merely giving it a descriptive label; the {\em year the big tree in the village blew down}; the {\em day the sun rises over that mountain}, etc. Enumerating rather than just labelling time permits one to know how long ago in the past something occurred, or how far into the future it will occur.
The Maya used their calendar to record astronomical events for astrological purposes \cite{thompson,saturno} and there continues to be much interest in understanding the Maya concepts of time and on what astronomical observations it may have been based \cite{chanier}.
It has been asserted that the Maya numeration system would be superior to today's, at least for the ease of recognizing small divisors of large numbers \cite{bietenholz}.
They certainly inherited parts of their number system, such as base 20, which was common across Mesoamerica, from earlier civilizations such as the Olmecs, 
and shared these aspects with succeeding peoples of Mesoamerica such as the Aztecs \cite{harvey,williams}.

In Mesoamerica there emerged a concept of zero, at first as a placeholder (an ordinal zero), before entering into arithmetic (a cardinal zero)  \cite{blume,mcneill}. As we shall discuss below, it is questioned whether the concept of zero was another such inheritance from the Olmecs to the Maya \cite{justeson}.  It is possible that the concept of zero has been discovered only twice: once in the Old World, where it  seems to have first appeared as a placeholder in Sumerian Mesopotamia four to five thousand years ago,  and once in the New. What seems certain is that the New World discovered zero on its own and that it was the Maya who fully developed the idea into a cardinal zero, used for calculations. 

A significant characteristic of the Maya calendar is the concurrent use of three separate number systems: the Haab, the Tzolk'in and the Long Count; the former two formed the Calendar Round, in which all dates are repeated every 52 years \cite{orozco,chanier}. 
 This combination of calendars is similar to our use today, without a second thought, of a year-month-day calendrical system together with an incommensurate week system, where repetition comes after 28 years\footnote{Disregarding the complications introduced in the Gregorian Calendar where centuries are not leap years unless divisible by 400.}.  The Maya civil calendar, the Haab, represents an annual solar cycle of 365 days, composed of 18 months (winals) of 20 days (kins) each, plus --- as in many calendar systems --- five extra epagomenal days  at the end of the year,  which were called unlucky days or days without name (wayeb, for the Maya) \cite{broda}.
 On the other hand, the divine calendar, the Tzolk'in, used to determine the time of religious and ceremonial events and for divination, has a cycle of 260 days, composed of 20 periods of 13 days each (trecena)  \cite{richards}, possibly owing to the 260-day span of time between zenithal sun positions at the latitude of  15$^{o}$N in Mesoamerica \cite{malstrom}.
  
Much debate has focused on who developed these calendrical systems. 
The Isthmus of Tehuantepec has long been seen as an important area of elaboration and differentiation of the first calendar in Mesoamerica, although opinions differ as to which side of the isthmus can claim precedence. Some scholars look to the Olmec society on the north side of the isthmus on the coast of the Gulf of Mexico, in the modern Mexican states of Veracruz and Tabasco. Others look to the south, to the Pacific coast, present-day Chiapas (Mexico) and Guatemala. And others indicate west, to modern-day Oaxaca \cite{rice}. 
Grove \cite{grove69,grove70} points out that the numeral glyph found in the Olmec culture, with a 260-day count in its calendrical inscription, may be the oldest. Similarly, Edmonson \cite{edmonson} proposes as the oldest calendrical record, one corresponding to the Olmec culture in the year 679 BCE.  Diehl \cite{diehl} indicates that in the decadence of Olmec culture, the epi-Olmec period, in Chiapa de Corzo in Chiapas and Tres Zapotes in Veracruz stelae were erected with the earliest known inscriptions of the Long Count.  Blume \cite{blume} notes that the earliest Mesoamerican Long Count is inscribed on Stela 2 at Chiapa de Corzo with a date of 7.16.3.2.13, corresponding to 36 BCE.\footnote{The date corresponding to the beginning of the Long Count, 0.0.0.0.0,  is Monday, 11th August, 3114 BCE,  according to the most accepted Goodman--Mart\'{\i}nez--Thompson correlation with our Gregorian calendar. It is supposed that this initial value was decided a posteriori in a similar fashion to how the current calendar era was proposed by  Dionysius Exiguus in the sixth century and widely implemented by Charlemagne in the 9th century. Similarly,  we have Unix epoch time, the number of seconds (minus leap seconds) elapsed since  00:00:00 UTC on 1st January 1970, retrospectively decided upon in 1972. Compare the earliest known numerical calendar, the Seleucid Era, which did begin with year 1 in 312/11 BCE \cite{kosmin}.}
In terms of  the development of mathematical ideas, we may affirm that the epi-Olmec and the proto-Maya  came together something over 2\,000 years ago in this fertile crescent and the Long Count, and zero, were the eventual results.
  
For everyday activities, the Maya used a pure vigesimal, base-20, numeral system (although there are no extant Maya documents showing this, and we know it only from what bishop Diego de Landa told of cacao bean counting in sixteenth century Yucat\'an)\footnote{\emph{``Que su cuenta es de V en V, hasta XX, y de XX en 
XX hasta C, y de C en C hasta 400, y de CCCC en CCCC 
hasta VIII mil. Y desta cuenta se serv\'{\i}an mucho para 
la contrataci\'on del cacao. Tienen otras cuentas muy 
largas, y que las protienden \emph{in infinitum}, cont\'andolas 
VIII mil XX vezes que son C y LX mil, y tornando \'a 
XX duplican estas ciento y LX mil, y despu\'es yrlo ass\'{\i} 
XX duplicando hasta que hazen un incontable n\'umero: cuentan en el suelo \'o cosa llana''} \cite{landa}. Landa's account, written circa 1566, in which he does not mention zero, demonstrates that Hindu-Arabic numerals were still being used little in Europe. Landa, who ordered the burning of almost all Maya texts, probably did not appreciate that he was destroying one zero in America just as another was struggling to emerge in Europe from the hegemony of Roman numerals. Ironically,  Landa's original manuscript is also lost, so we cannot be sure that the version we have with Roman numerals is how he wrote it.} \cite{landa}.
In their Long Count calendar (Figure~\ref{fig:Maya}), however, they would use a slightly modified version of this. The first and second place values were $20^0$ and $20^1$ as usual, but the third was $20 \times 18$. This is presumably because $20 \times 18 = 360$ represents much more accurately than $20^2 = 400$ the number of days in a year. All subsequent place values were multiplied by 20. Thus we have 1 kin (= 1 day), 1 winal = 20 kins, 1 tun = 18 winals, 1 katun = 20 tuns and 1 baktun = 20 katuns.
Accordingly, a number would be expressed in this system as 
$$ N = d_k (18 \times 20^{k-1}) +\ldots + d_3 (18 \times 20^2) + d_2 (18 \times 20) + d_1  20 + d_0.  $$
This is an example of a non-power positional number representation system\footnote{Chrisomalis  \cite{chrisomalis2010} has argued
that \emph{``there was no Maya positional numerical notation system''} because extant evidence shows the Maya using the Long Count only for measuring time, and because --- he maintains --- \emph{`\,``9 millennia, 4 centuries, 3 decades, 6 years'' is read and understood differently from ``9436 years'''}, and he opines that the former, unlike the latter, is not a positional number representation system. Whether or not it was the case that the Maya only used the Long Count system in this one setting of the calendar (the absence of evidence may be owing to the destruction of most Maya codices), however, we disagree with his latter assertion, because the key point here is that the Long Count possesses the formal structure of a positional number representation system. There are manifold instances in the history of mathematics in which the general utility of a formal structure has not become apparent to mathematicians until long after its  invention. This may, or may not, be another one of those cases.}. 
Since the digits go up beyond 9, to avoid using extra non-decimal digit symbols to write Long Count numbers the convention is to use the following notation with intercalated dots between digits written in decimal
to avoid any confusion between numbers:
$$ N = d_k.\cdots.d_1.d_0.$$

\begin{figure}
 \centering a) \includegraphics[height=0.5 \columnwidth]{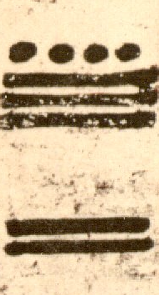}\,\,b) \includegraphics[height=0.5 \columnwidth]{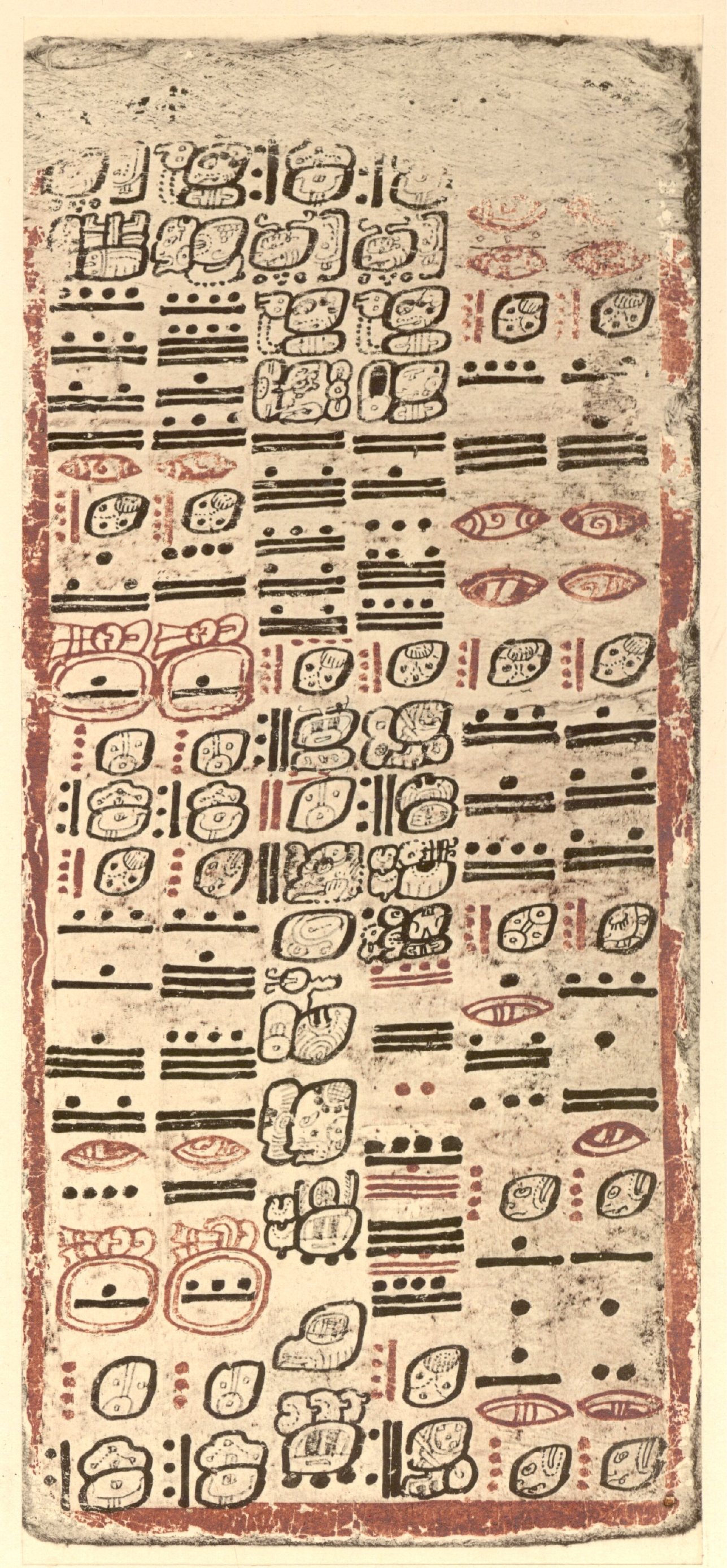}
\caption{
Instances from the Dresden codex of Maya numbers written with  an 18 or a 19 in the second place.
The dot and bar notation of the Maya seen here composes the digits 1--19 using zero to four dots to represent ones and zero to three bars to represent fives, so that 1 is one dot, and 19 is three bars below four dots. Generally numbers were written vertically with the most significant figure at the top.
(a) 390 written as 19.10 rather than as 1.1.0 on page 72;
(b) 10.11.3.19.14 or 10.11.3.18.14, i.e., 1\,520\,654 or 1\,520\,674 in decimal, on page 70. It is unclear whether is there a dot missing in the second place digit owing to wear of the codex, but the unequal dot spacing --- compare with 19 as written in (a) --- makes it plausible that there were originally four dots, with this digit thus reading 19 rather than 18.
\label{fig:dresden}      
}
\end{figure} 

In a regular base-20 system, when a place value is completely filled, we simply write 0 and carry a 1 to the next power. 
For instance, we can  fill up the units place with numbers 1 up to 19, but on reaching 20 we have to write it as $1.0 =  20^1 + 0 \times 20^0$, and similarly for higher powers.
In the calendar count, however, the third place being $18 \times 20$ creates some difficulties. If the second place ($20^1$) is filled up, we would have 20 sets of 20 which cannot be carried over to the next power: only $18 \times 20$ can, leaving $2 \times 20$ back. The same happens if we have $19 \times 20$ in the 20's place: it can be carried over leaving $1 \times 20$ behind.  
So for example, $3.19.3 = 4.1.3$ and $7.18.11 = 8.0.11$. Thus, if the digits can go up to 19 in the second place, this non-power positional number representation would not be unique starting from Maya Long Count numbers from \(18.0=1.0.0\) to \(19.19=1.1.19,\) corresponding to the base-10 numbers 360  to 399. In other words, the number system is partially nonunique, with the nonuniqueness affecting about 10\% of numbers.

So did the Maya ensure uniqueness in the Long Count by having the second digit go only to 17? That might make sense considering that this number representation system was used as a calendar (recall 18 20-day `months'  plus five extra days made up the Maya year).  In that context, one might naturally carry directly into the larger units. 
Indeed, Freitas and Shell-Gellasch  \cite{freitas} did not find any examples: 
``no Maya Long Count numbers with an \(18\) or \(19\) in the second place appear on known monuments or documents'' 
they wrote. 
But this is not so. 
Closs \cite{closs} notes an example in the Dresden codex, of $390_{10}$ (the subscript refers to the base) expressed in the form 19.10.  We show this instance in Figure~\ref{fig:dresden}a.  Note that Closs was expecting the number  to appear as 1.1.10, and was surprised to find it written in this other manner.  Cauty and Hoppan \cite{cauty} found this same example and noted a further instance in the Dresden codex, which we show in Figure~\ref{fig:dresden}b, where instead of 10.11.4.1.14 there is written 10.11.3.19.14. 
Again, they see these instances as being irregular variants where the scribe has omitted to carry into the larger units.

\begin{figure}
  \centering\includegraphics[width=0.7\columnwidth]{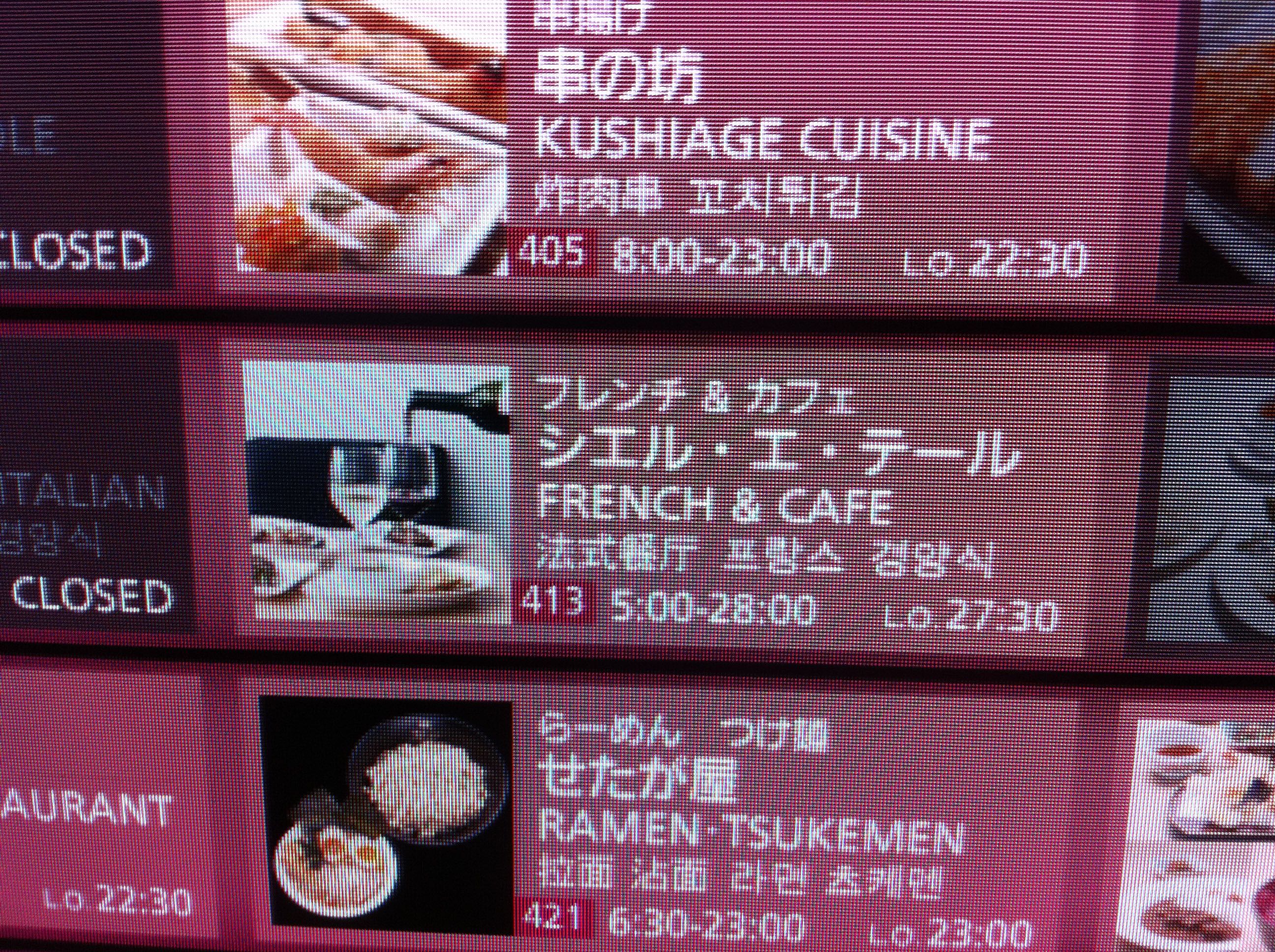}
\caption{Today one occasionally finds an example of Maya-style notational non-uniqueness in our representation of time.
Information on restaurant opening hours at Tokyo Haneda airport; notice the opening hours and time for Last Orders of the second entry.
\label{fig:airport}      
}
\end{figure}

We may note that many old units of money and measure functioned in this same way.
E.g., in the old British monetary system (imported from Charlemagne's continental Frankish empire) of pounds, shillings and pence, 12 pennies made a shilling and there were 20 shillings in a pound:
$ N =   d_2 (20 \times 12) + d_1  12 + d_0$,
where the digits $d_0$, $d_1$ go only to 11, 19; that is, to $b_{i+1}/b_i-1$. Likewise  the linear measures 
12 inches = 1 foot, 3 feet = 1 yard, 220 yards = 1 furlong, 8 furlongs = 1 mile etc, ensured uniqueness by having the digit go only to 1 less than the next level \cite{fraenkel}. 
And of course we do precisely the same with our modern calendar and our timekeeping: we write dates as days ($<$ a month), months ($<$ a year) and years, and the same with hours, minutes ($<$ an hour) and seconds ($<$ a minute). We would not generally think of giving a date as 13/13/2018, rather than 13/1/2019, nor a time as 12:65 rather than 13:05; but see Figure~\ref{fig:airport}. However, one instance where we are perfectly comfortable with such nonuniqueness today is in currency, where coins and notes in denominations often based on 1, 2, 5 permit us to pay a given sum of cash in multiple ways.

It is notable that we are worse off today with our calendar with a jumble of months with different lengths than the ancient Maya with their 20-day months, as we have to remember that ``30 days hath September...'' etc, in order to perform calendar calculations with our very irregular length months inherited from the Romans. Moreover, although we are happy to consider the first minute of the hour, minute zero, and the last minute, 59, only in the 24-hour clock do we condescend to have an hour zero, and we refuse to consider a day zero or a month zero, just as there is no year zero in today's Gregorian calendar. 
It should also be noted that the Tzolk'in --- perhaps the earliest Mesoamerican calendar \cite{blume}  --- has days 1 to 13, without a zero in the same way as the days of our months lack a zero day.  On the other hand, the days of the solar Haab calendar use the same digit notation of dots and bars as we see in Figure~\ref{fig:dresden} used for coefficients of Long Count quantities, so it is natural to ask how the Maya got to be more logical than us, to arrive at a day zero.

\section{Bijective Numeration}

It is not widely appreciated that a positional number representation system does not need a zero.  Instead of digits in the range from $0$ to $b-1$, we can simply shift them by one to the range from $1$ to $b$ \cite{foster} (see the Methods section for a discussion of digit shifting).  We have done away with zero --- whose introduction is often held to have been essential for the development of positional number systems --- yet we can still represent all numbers uniquely. For instance, if we do this with base 10, we simply need a new digit symbol for ten; let us borrow from the Romans and use X. Then the digits are from 1 to X, and most numbers remain written as in normal base 10. Only those containing 0 are altered: 10 becomes X, twenty, 1X, one hundred, 9X, and so on (Figure~\ref{zero_fig}a). It is still the case that $1 + 1 = 2$, but now $9+1 = X$ \cite{knuth}. This zero-less number system has sometimes been called bijective numeration \cite{allouche}; it is bijective or one to one because in this system there is no possibility to have leading zeros in front of a number. It keeps on being rediscovered \cite{forslund,boute,manca}\footnote{
Note that one may combine bijective numeration with non-power number representation systems in the same way as with a fixed base system (Figure~\ref{zero_fig}c).} 
We see, then, that the introduction of zero, although viewed historically as linked to the development of our Hindu-Arabic decimal positional number system, was not necessary for a positional number system.  

\begin{figure}
\centering\footnotesize
a)$\begin{array}{c|c}
\textrm{Bijective base 10}  & \textrm{Decimal} \\
\hline
1 & 1 \\
2 & 2 \\
$\vdots$ &$\vdots$ \\
9 & 9 \\
X & 10 \\
11 & 11 \\
12 & 12 \\
$\vdots$ &$\vdots$ \\
19 & 19 \\
1X & 20 \\
21 & 21 \\
$\vdots$ &$\vdots$ \\
99 & 99 \\
9X & 100 \\
X1 & 101 \\
X2 & 102 \\
$\vdots$ &$\vdots$ \\
XX &110 \\
111 & 111 \\
$\vdots$ &$\vdots$ \\
\end{array}$ 
b)$\begin{array}{c|c|c}
\textrm{Bijective base 20}  & \textrm{Base 20 with zero} & \textrm{Decimal} \\
\hline
1 & 1 & 1\\
2 & 2 & 2 \\
$\vdots$ &$\vdots$ &$\vdots$\\
(19) & (19) & 19 \\
(20) & 10 & 20 \\
11 & 11 & 21 \\
$\vdots$ &$\vdots$ &$\vdots$ \\
(19)(19) & (19)(19) & 399 \\
(19)(20) & 100 & 400 \\
(20)1 & 101 & 401 \\
$\vdots$ &$\vdots$ &$\vdots$ \\
(20)(20) & 110 & 420 \\
111 & 111 & 421 \\
$\vdots$ &$\vdots$ &$\vdots$ \\
\end{array}$ \\
\vspace{0.5cm}
c)$\begin{array}{c|c|c|c}
\textrm{Bijective Long Count} & \textrm{Long Count with zero \& twenty} & \textrm{Long Count with zero}  & \textrm{Decimal} \\
\hline
1 & 1& 1 & 1 \\
2 & 2 & 2 & 2 \\
$\vdots$ &$\vdots$ & $\vdots$ &  $\vdots$\\
19 & 19 & 19 & 19 \\
20 & 1.0 \textrm{ or } 20 & 1.0 & 20 \\
1.1 &  1.1 & 1.1 & 21 \\
$\vdots$ &$\vdots$ & $\vdots$ &  $\vdots$\\
17.19 & 17.19 & 17.19 & 359 \\
17.20 & 1.0.0 \textrm{ or } 17.20 \textrm{ or } 18.0 & 1.0.0 \textrm{ or } 18.0 & 360 \\
18.1 & 1.0.1 \textrm{ or } 18.1 & 1.0.1 \textrm{ or } 18.1 & 361\\
18.2 & 1.0.2 \textrm{ or } 18.2  & 1.0.2 \textrm{ or } 18.2 & 362 \\
$\vdots$ &$\vdots$ & $\vdots$ &  $\vdots$\\
18.20 &  1.1.0  \textrm{ or }  18.20  \textrm{ or } 19.0 \textrm{ or } 1.0.20 & 1.1.0  \textrm{ or } 19.0  & 380 \\
1.1.1  \textrm{ or }   19.1 & 1.1.1  \textrm{ or }   19.1 & 1.1.1  \textrm{ or }   19.1  & 381 \\
$\vdots$ &$\vdots$ & $\vdots$ &  $\vdots$\\
1.1.10  \textrm{ or }   19.10 & 1.1.10  \textrm{ or }   19.10 & 1.1.10  \textrm{ or }   19.10  & 390 \\
$\vdots$ &$\vdots$ & $\vdots$ &  $\vdots$\\
1.1.20 \textrm{ or } 19.20   &  1.2.0 \textrm{ or } 19.20 \textrm{ or } 20.0 \textrm{ or } 1.1.20 & 1.2.0 & 400 \\
1.2.1 \textrm{ or } 20.1  &  1.2.1 \textrm{ or } 20.1 & 1.2.1 & 401 \\
$\vdots$ &$\vdots$ & $\vdots$ &  $\vdots$\\
1.2.20 \textrm{ or } 20.20  & 1.3.0  \textrm{ or } 20.20  \textrm{ or } 1.2.20  & 1.3.0 & 420 \\
1.3.1 & 1.3.1 & 1.3.1 & 421 \\
$\vdots$ &$\vdots$ & $\vdots$ &  $\vdots$\\
 \end{array}$
\end{figure} 
\begin{figure}[t] 
\caption{
(a) Bijective base 10 with digits 1--X  differs from base 10 with a zero --- our usual decimal system with digits 0--9 --- only when we write the numbers that would have a zero in decimal; otherwise the two are the same.
(b) Similarly, to move between a bijective numeration system with digits 1--(20) and one using a zero, with digits 0--(19) (which of course would require symbols for the digits (10)--(19) or (20), which are written for this reason with parentheses around the number), the differences are only with numbers involving a zero.
The conversion works similarly with any base, and 
(c) with mixed base systems, as we see with the different versions of the Long Count. The bijective Long Count with digits 1--20, the Long Count with digits 0--20, and a version with digits 0--19, can easily coexist. Most numbers are written the same in all three versions. Only those numbers involving digits 0 and 20 in the first place, and 0 and 18 onwards in subsequent places, can differ. There is redundancy in the three versions, the same but shifted in the 0--19 and 1--20 versions, and the combination of these two sets of redundancies in the 0--20 version.
}
\label{zero_fig}      
\end{figure}

A key point in exploring the Maya number systems is this: given a positional number system that uses an unfamiliar set  of symbols, how may we know whether the symbols include a zero, or not?  
How do we know whether the Maya in their Long Count were writing their calendar using days from 0 to 19, or from 1 to 20? Clearly this question is pertinent given that today our calendar does not include a day zero. That is to say, did the Maya really begin their months with day zero and end them with day nineteen, or did they begin with one, like us, and end with twenty?  If the meaning of the digit symbols is completely unknown a priori then only way to answer this question is to look at arithmetical operations with these symbols \cite{forslund}\footnote{It reminds us of the passage in \emph{Alice's Adventures in Wonderland} \emph{``I'll try if I know all the things I used to know. Let me see: four times five is twelve, and four times six is thirteen, and four times seven is  --- oh dear! I shall never get to twenty at that rate!''} \cite{carroll1865}, which can be understood if Alice is counting in a varying base: she is expressing $4n$ in base $3n+3$, and she cannot get to 20 --- as many people have pointed out --- since after $4\times 12 = 19_{39}$, $4\times 13 = 1X_{42}$.}. 
In the four surviving Maya codices, written on astronomical and calendrical themes between the classical Maya period and the arrival of the Europeans, a shell symbol --- thought to be a stylized image of the shell of a gastropod mollusc of the genus \emph{Oliva} \cite{blume} --- represents zero. That it does represent zero is clear from, for example, multiplication tables in the Dresden codex that would be incorrect if we tried to interpret the shell symbol as 20. (Although it should be noted that the arithmetic in the codices contains errors.) So, certainly by their post-classical period during which the codices were written, we do have a zero-based positional system, with digits 0--19. 

There is a variety of evidence pointing towards an earlier bijective system, with digits 1--20, both within the codices themselves and in the stone inscriptions on the earlier classical stelae.
Within the codices there is evidence for an earlier, non-positional system with an explicit symbol for twenty. In the Dresden Codex, distance numbers between 20 and 39 are frequently expressed by prefixing with the dot and bar notation a number between 0 and 19 to a moon glyph representing twenty, and similar ``x+20'' notation was used earlier on stelae, as has been discussed by Thompson \cite{thompson1950}.
Although an explicit zero is first used for a Long Count inscription of 8.16.0.0.0 (357 CE), at Uaxactun on Stelae 18 and 19, a later 711 CE inscription at Stela 5 at Pixoy is not written as 9.14.0.0.0, but as 9.13.20.0.0, as Closs has pointed out \cite{closs1978}.
And at the Temple of the Cross (Palenque, Chiapas) there are the forms 20 Mol (in 13 Ik 20 Mol = 13 Ik 0 Ch'en) and 0 Zac (in 9 Ik 0 Zac) entered in two Haab dates written side by side  \cite{cauty}.

\section{The Development of Zero}

The route to zero that Mesoamerica took must be teased out from the sparse evidence available. Mathematics can help with this task. Initially there was a non-positional number system with digits from 1--13, without a zero, in the Tzolk'in calendar. Archeological evidence of this calendar has been found in Olmec cave paintings dated 800--500 BCE \cite{grove69,grove70}.

Then there is the Haab, again a non-positional  number system. The Haab may possibly have been set up around 500 BCE \cite{bricker1982}; there is archeological evidence from 500--400 BCE from Monte Alb\'an,  Oaxaca \cite{marcus1976}.  As we have indicated, the days of the Haab generally run not from 1--20, but from an initial day, followed by the 1st, 2nd, etc, up to the 19th day. And as we have pointed out above, in the Haab on occasion a glyph for the end of a month, i.e., for day 20, was used instead of that for the beginning of the month. 
Of course then  it is natural today to translate the glyph for the initial day, often referred to as \emph{chum}, by zero. However, it is not one of the same glyphs as the zero later found in the Long Count. Maya scholars have debated for many years about the meaning of chum for the Maya \cite{blume}. Some have thought it to be the end of the preceding month, i.e., a species of 20, and others the beginning of the new month, i.e., a type of zero.  Some translate it in a non-numerical way as the ``seating of'' the month \cite{bowditch1901,gordon1902,thompson1950}.
What we can put forward for our purposes about chum with relative certainty, however, is twofold: (1) that it does not perform the same function as the zero in the Long Count, being an ordinal, not a cardinal zero; 
and, however, 
(2) that it may well have influenced the development  of  the zero placeholder to come in the Long Count. 

Now let us move to the Long Count, for which the earliest evidence is several hundred years later than the two preceding calendars. 
The eight earliest known Long Count inscriptions (which are generally designated epi-Olmec rather than Maya \cite{diehl}) are 7.16.3.2.13 (Chiapa de Corzo Stela 2), 7.16.6.16.18 (Tres Zapotes Stela C), 7.18.9.7.12 or 7.19.15.7.12 (El Ba\'ul Stela 1), 8.3.2.10.15 and 8.4.5.17.11 (Takalik Abaj Stela 5)\footnote{Davletshin \cite{davletshin2002}, Justeson \cite{justeson2001}, and Macri \cite{macri2020} have proposed a different interpretation of the dates inscribed on Takalik Abaj Stela 5. In their interpretation, the inscriptions, which read as 8--3--2--10--5--[damaged day sign] and 8--4--5--17--11--[damaged day sign] would correspond to Long Counts 8.3.2.[0].10 and 8.4.5.[0].17, where [0] represents an implicit zero. If this alternative interpretation is correct, an implicit zero here also supports our arguments.}, 8.5.3.3.5 and 8.5.16.9.7 (La Mojarra Stela 1) and 8.6.2.4.17 (the Tuxtla statuette) \cite{malmstrom1996}. It is tantalizing that  in none of the eight known inscriptions that show the earliest development of positional notation in the Long Count  can we find the full range of digits that were necessarily then in use, which would enable us to understand whether the calculations behind the Long Count were being performed with digits 1--20, 0--19, or 0--20. Neither 0 nor 20 appears in any of the earliest examples written from 36 BCE to 162 CE. We know that there necessarily has to have been one of these three systems in use in order to satisfy the mathematics of the Long Count; in order to be able to write all numbers (Figure~\ref{zero_fig}c). With fewer digits with the same multipliers not all numbers can be represented and so some calculations simply could not be performed.

Then there appear in the archeological record examples with an implicit zero denoted by the lack of a digit, before finally we get examples with the explicit written cardinal zero, as well as the example with an explicit twenty on Stela 5 at Pixoy. (Since the  Tzolk'in and Haab continued to be used alongside the Long Count, we can use this calendrical redundancy for error checking, to make sure that we understand correctly whether a glyph is a zero or a twenty.) As in the case of the solar calendar glyph  \emph{chum}, within Maya scholarship there has been a great deal of debate about how to read the Long Count glyphs that are not bar and dot numerals. In contradistinction to the case of the solar calendar, here we can check Maya arithmetic and see that we really do have a glyph for zero in the two forms in common use, one a species of quatrefoil, and the other a hand. By the time we have archeological evidence these are glyphs for zero, but might these same glyphs once have been glyphs for twenty, whose meaning had altered over time? That is an idea that has occurred to a number of Maya scholars over many decades \cite{bowditch1901,gordon1902,thompson1950,blume}. 

In terms of the mathematics, when altering a digit from a twenty to a zero, one is moving from bijective numeration to a non-bijective system which merely requires incrementing the digit in the superior position by one, so the idea of a shift in meaning of these glyphs is quite tenable.
Moreover,  in fact one can use the system with digits 0--20, with both a 20 and a 0, without any confusion, as we show in  Figure~\ref{zero_fig}c. 
Given that at its first sightings the zero is present only implicitly, as an absence --- as a missing digit --- we can understand that it was a new concept whose usage caused deep conceptual problems,  in a similar way to how irrational numbers, imaginary numbers, non-Euclidian geometries, etc, have caused problems at various times to mathematicians. All of which strengthens the idea of a previous usage of an explicit twenty in the system. Although --- barring fresh archeological finds --- we do not see the first part of the process with digits 1--20 in extant inscriptions, we see precisely this latter stage with digits 0--20 in the Long Count written on Stela 5 at Pixoy.

Thus we see the mathematics of the development of the Long Count as being one in which a non-power positional system was being derived, with the mixed multipliers of 20 and 18 that we have described --- possibly alongside a pure base-20 system like Figure~\ref{zero_fig}b for civil use --- together with a set of digits that initially ran from 1--20, and then, much more nontrivially, was changed to include the possibility of absence: of a zero\footnote{It is notable that Justeson \cite{justeson}, who comes at this question from a completely different approach to our own, arrives at a similar conclusion on this point.}. All these digits and multipliers could be used together because the Maya mathematicians were happy with the flexibility of their number representations leading to non-uniqueness, to redundancy, so that they could write both a zero and a twenty side by side in a number such as 9.13.20.0.0 at Pixoy.  That is to say, their positional number system broke both sufficient conditions for uniqueness listed in the Methods section: neither are the positional weights powers of a base $b$, nor are the digits limited to a range from 0 to $b-1$, or to a shift of that range.  The Long Count with zero and twenty gradually gave way to a Long Count with zero.  By the time  the extant codices were written, when Maya civilization was on the wane,  the Long Count was worked almost always with digits with the explicit zero, using the shell glyph which is characteristic of the codices, and the use of twenty as a digit had become vestigial.

\section{Ex Nihilo Nihil Fit?}

In the Old World, zero first emerged from the development in Sumerian Mesopotamia of a sexagesimal positional number system for accounting purposes. This was to begin with an implicit zero, by which we mean that at first it had no symbol associated with it, but was simply a lack of a digit. This makes perfect sense within the scheme of tallying goods: a lack of something corresponds to a missing number in its corresponding column. But it became awkward that the lack of something could be misinterpreted when writing tallies in columns, in a positional notation, and so after some time the implicit zero was given its own symbol and became explicit as a placeholder in the base-60 notation. So it was natural, from this Old World point of  view of counting goods, that a positional number system and a zero should go together. In the New World, however, the impulse for the positional number system came not from counting goods, but from the calendar, from counting days. And in counting days, what would be the lack of a day? It is much more natural to use the counting numbers, 1, 2, 3. We see this today in our Gregorian calendar  with no day, month, or year zero. So from the American point of view, it made sense that there could be a positional notation without a zero.

Mesoamerica did not have to discover zero at the same time as inventing positional notation; the two are independent concepts.
With a bijective positional number system one can represent all numbers, and one does not need a zero. One has to ask whether the Maya used it;  that is, was the symbol they used a 0 or a 20? 
That the Maya discovered the concept of zero beginning from a bijective non-power number system is a plausible hypothesis when one considers on one hand the available historical evidence, which has led Maya scholars to debate whether chum meant zero or twenty, and on the other hand the ease with which one can move from a bijective system to a system with a zero.
We infer that initially they used 20 and only later 0, and they shifted via an intermediate Long Count with both a twenty and a zero that we see in the historical record. To go between these different systems affects only the numbers whose representation contains the digit (20), in which one replaces a (20) with a (0), at the same time adding 1 to the digit in the superior position\footnote{Of course  
as the number gets larger, the probability of it containing a (20) tends to one. However, although on occasion the Maya did represent numbers with a large number of digits, most of the numbers written have five or fewer digits.}. If one adds to this mixture of mathematics and anthropology the point that the Maya, owing to the redundancies built in to their mixed-based system, were used to the idea of the same number being represented by more than one different sequences of symbols, one can understand that they could make the momentous conceptual leap to using a cardinal zero owing to their familiarity with multiple and redundant number representation systems. 

Zero is a slippery and difficult concept that took a long time to be accepted as a number with its own symbol to use in calculations. In Europe zero arrived in the Middle Ages as one of the Hindu-Arabic digits with the decimal positional number system, but these two ideas did not have to come at the same time. Zero is not needed for a positional number system. Different ways to write calendar dates developed in Mesoamerica two millennia ago helped to invent a number system in which the same number can be written in different forms even without a zero.  America discovered zero on its own and later it was added by the Maya to this, the Long Count.

\section{Methods}

The history of mathematics has an intrinsically interdisciplinary character. In order to make matters clear for a diverse readership, and to provide a reasonably self-contained argument, in this section we spend a little time in an outline of the mathematics involved in number representation systems.

\subsection{Number Representation Systems}

``The mirror of civilization'' is what Hogben termed mathematics in his \emph{Mathematics for the Million} \cite{hogben}, ``interlocking with man's common culture, his inventions, his economic arrangements, his religious beliefs''.
It may be that the initial use of the symbolic management of numbers through visual signs corresponded to utilitarian needs, for example, for the exchange of goods; the first form of commerce. However, numbers became part of the human endeavour for knowledge very early. Perhaps astronomy --- the counting of the elapsed time between recurring events of day--night, winter--summer, relative positions of planets and stars, eclipses, and so on --- was the earliest `scientific' application of number systems\footnote{Until the scientific revolution of the 15th--16th centuries, astronomers were also astrologers or priests, and astronomical data were used for astrological or religious rather than what we would think of as scientific purposes. Nonetheless, as numbers were used to document, explain and predict natural phenomena,  we may consider this a proto-scientific application.}. Geodesy also has represented an important practical aspect in early civilizations that led rapidly into geometric developments. We find the apex of number in the Pythagorean doctrine that the entire universe is governed by numbers; for the Pythagoreans, that meant integer and rational numbers \cite{cartwright2020}. 

The felicitous choice of a numeration system is relevant for solving specific problems and also for developing and improving mathematical models and algorithms \cite{bell,fraenkel}. We can see from our position of hindsight that  civilizations that used inconvenient number systems were held back in their development of mathematical knowledge. 
Today's decimal number system is a positional system, but there are many historical examples of non-positional numeration systems \cite{Menninger}. A familiar example is that of Roman numerals.
In the Roman system the values of the symbols are in general independent of their position; numerals are written from left to right in descending order, writing the biggest numeral possible at each stage. There is only a relative positional dependence that determines whether a particular number should be added to or subtracted from its neighbour for obtaining the represented number.  For example, I represents 1 and V, 5, and there are two ways to write 4, IV and IIII (with the latter version generally seen on clocks; the subtractive rules leading to forms like IV were alternatives that only became usual in later Roman periods).  Thus Roman numerals constitute primarily a non-place-value system, but because of the use of the subtractive principle --- e.g., IV represents four while VI represents six ---, the Roman system may be classified as a mixed system.
Arithmetical operations are very difficult to implement with non-positional systems. Cultures that used non-positional representation systems generally relied on mechanical aids for performing operations between numbers, the most popular of these being the abacus. On the contrary, positional number systems allow a compact representation and easy implementation of arithmetical operations. Mainly owing to this last feature, positional number systems have historically prevailed over non-positional ones. The current decimal representation system with Hindu-Arabic symbols representing the digits was spread in the West by Leonardo Pisano, better known as Fibonacci, with his 1202 book \emph{Liber Abaci} (the Book of Calculations) \cite{fibonacci}, and gradually replaced the cumbersome Roman system. In Mesoamerica a positional number system was in use much earlier.
 
Most of the numeration systems that we know and use are univocal, that is, they are not redundant; univocity means that any symbol represents only one number and, conversely, any number is represented by only one symbol (except for the minor point of leading zeros that we shall discuss below). This is usually achieved using a power positional system, defined by a small set of integers given symbols, called digits, that, depending on their position along a representation string, are multiplied by the powers of a given base or radix. However, there have existed in the past, and there continue to exist, non-power positional number representation systems in which the multipliers are not the powers of a given base. Such mixed-radix systems were studied by Cantor \cite{cantor1869}.
One historical example in which the necessity to describe more adequately annual timescales led to a non-power number representation system is the Maya Long Count calendar, depicted in Fig.~\ref{fig:Maya}.  Others are the old systems of money, weights and measures from around the world not based on multiples of ten. While these systems were often employed in such a way as to preserve univocity, non-power number representation systems, mainly employed in computing for transmission and storage of data, are today used so as to be redundant.  The reason for this redundancy is to have error detection, and so the possibility for error correction, built in to the system \cite{diego2008}\footnote{The redundancy of overspecifying calendrical information performs precisely this function of error detection in our present calendar, as for example when we write Saturday 17th November 2018.  (Algorithms to determine the day of the week for any given date have been devised by mathematicians from Gauss \cite{gauss}, through Lewis Carroll \cite{carroll}, to John Horton Conway \cite{conway}.) We can presume that Maya scribes understood and used this same redundancy for error checking purposes when they wrote dates using the three Maya calendars. The same date in the Long Count is 13.0.5.17.17, 3 kab'an, 10 Keh.}. It is thus most interesting that in DNA, the biological molecule of information  transmission and storage, we find in the genetic code the structure of a non-power number representation system \cite{gonzalez2004,gonzalez2016}.

\subsection{Positional Number Systems}

A counting system is said to be {\it positional} if each digit is weighted with a different value according to its location in the string.  
The most common positional numeral systems are power representation systems where the positional weights are powers of some number $b$, called the base or radix, and the digits are allowed to take any value from $0$ to $b-1$.
The main advantage of such a system is that any integer $N$ has a unique representation of the form
$$N = d_k b^k + d_{k-1}b^{k-1} + \ldots + d_0 b^0,  \qquad  0 \leq d_i \leq b-1 \quad \forall i  . $$
The decimal system, base 10, is undoubtedly the most familiar and widespread example, but it is not the only one. The first place-value system, developed by the Mesopotamians, was  sexagesimal, base 60, which is why we still measure angles and time in units of 60. In more recent times, the binary, base 2, system has become a fundamental tool in informatics; base 16, hexadecimal, and base 8, octal, are also important in computing.

\subsubsection{Uniqueness of the representation}

As previously stated, in a positional system with base $b$ each number has a unique representation.
This can be proved by contradiction, i.e., assuming that an integer $N$ can be written in two different ways: 
\begin{align}
 N &=  d_0 b^0 +  d_1 b^1 +\ldots +  d_k b^k, \qquad   0 \leq d_i \leq b-1\\
 N &=   a_0 b^0 +  a_1 b^1 +\ldots +  a_k b^k, \qquad  0 \leq a_i \leq b-1
 \end{align}
and assuming that $\exists i$ such that $d_i \neq a_i$. In particular, assume $a_i > d_i$, and $\forall j > i, d_i = a_i $.
Subtracting (2) from (1), 
 $$0 = (d_0 - a_0) b^0 + (d_1 - a_1) b^1 + \ldots + (d_{i} - a_ {i})b^{i} $$
$$\Rightarrow  (a_i - d_i) b^i = (d_{i-1} - a_ {i-1})b^{i-1} + \ldots + (d_0 - a_0) b^0 .$$
By assumption, $a_i - d_i > 0 $ and so
\begin{equation*} 
b^i \leq (a_i - d_i)b^i  \Rightarrow b^i \leq (d_{i-1} - a_ {i-1})b^{i-1} + \ldots + (d_0 - a_0) b^0 .  
\end{equation*}
Since $d_j$ is a digit, $d_j \leq b-1$ and thus  $d_j - a_j \leq b-1 \quad \forall j$.
So one obtains
\begin{equation*}
\begin{aligned}
  b^i &\leq (d_{i-1} - a_ {i-1})b^{i-1} + \ldots + (d_0 - a_0) b^0 \\ 
  &\leq (b-1)b^{i-1} + \ldots + (b-1) b^0 
 \end{aligned}
 \end{equation*}
$$\Rightarrow b^i \leq (b-1)(b^{i-1}+ \ldots + b^0) .$$
But $(b^{i-1}+ \ldots + b^0)$ is a geometric series, and we know that
$$ \sum_{j=0} ^{i-1} b^i = \frac{b^i -1}{b-1} ,$$
and so we obtain
$$  b^i  \leq b^i-1,$$
 which is a contradiction;  the proof is complete. 

\subsection{Signed-Digit Representation}

It is clear from the above proof that the uniqueness of the power representation is given by two key assumptions:
\begin{enumerate}
\item the digits are limited to a range varying from $0$ to $b-1$, and
\item the positional weights are the powers of $b$.
\end{enumerate}
We can  develop non-standard positional numeral systems, where one of these two conditions is not fulfilled.

A first relaxation of the above conditions is where digits are allowed to go beyond the prescribed range. A particular example is the signed-digit representation, where each digit is given a positive or negative sign, hence its name. Uniqueness cannot be guaranteed anymore; it is easy to see that this representation is in fact redundant. Let us take for instance the binary signed-digit case: the positional weights are powers of 2, just like the usual binary, but the digits can now take values $-1, 0$ and $1$. The number 9, for example, can be written in three different ways: 
$$ 1001 = 11\bar{1}\bar{1} = 101\bar{1} $$
with the convention that $\bar{1} = -1 $.

Redundancy can be eliminated by considering the so called {\it balanced form} of the representation. Given a base $b$, the allowed digits are $b-1$ numbers $d_{i}$ taken from the range
$$  - \left \lfloor{\frac{b}{2}}\right \rfloor \leq d_{i} \leq (b-1) - \left \lfloor{\frac{b}{2}}\right \rfloor $$
where the floor function $\left \lfloor{x}\right \rfloor$ maps $x$ to the largest integer smaller or equal to it.
For simplicity, we will only consider the case where $b$ is an odd integer, which implies
$$\left \lfloor{\frac{b}{2}}\right \rfloor = \frac{b-1}{2} $$
and hence $$  -\frac{b-1}{2} \leq d_{i} \leq \frac{b-1}{2} .$$
We can prove that the balanced form is unique in the following way. 
We have shown above that every integer has a unique representation in base $b$ of the form
 \begin{equation} N = d_n b^n + d_{n-1}b^{n-1} + \ldots  + d_0 b^0,  \qquad  0 \leq d_i \leq b-1 \quad \forall i . \end{equation}
First consider the coefficients $d_{k} = b-1$.
Noting that
$$d_{k} b^k = (b-1)b^k =  b^{k+1} - b^k $$
one can substitute this expression for $d_{k} b^k$ into (3) to get
\begin{align*} 
&N = \\
& d_n b^n + d_{n-1}b^{n-1} + \ldots +d_{k+1}b^{k+1}+b^{k+1}+ (-1)b^k+ \ldots + d_0 b^0\\
&= \\
& d_n b^n + d_{n-1}b^{n-1} + \ldots +(d_{k+1} + 1)b^{k+1}+ (-1)b^k+ \ldots  + d_0 b^0 .
\end{align*}
If $d_{k+1} + 1 = b-1$, we repeat the previous step until there are no more coefficients equal to $b-1$.
Then, we seek to eliminate every coefficient of the form $d_t = b-2$, replacing it with $$ (b-2)b^t = b^{t+1} - 2b^t .$$
Finally, we get to the digits of the form 
$$ \left (b- \frac{b-1}{2} \right) b^s = b^{s+1} -\left(\frac{b-1}{2}\right)b^s $$
and we plug this expression into $N$.
Therefore, with this process we have found a unique representation for $N$ with digits drawn from the range 
 $$  -\frac{b-1}{2} \leq d_{i} \leq \frac{b-1}{2} $$
as required. 

In intuitive terms, what we have done is simply to shift down the allowed interval for the digits; the proof shows that in this case uniqueness is preserved.
For example, consider a normal base-3 representation (observe that the allowed digits are 0,1,2).
Take the number $$211_3 = 3^0 + 3^1 + 2 \times 3^2,$$ where the subscript refers to the base.
In order to reduce it to a balanced form representation, we seek to have only $-1,0,1$ as digits. We eliminate the digit 2 as follows:
 \begin{equation*}
\begin{aligned}
   3^0 + 3^1 + 2 \times 3^2 &=3^0 + 3^1+ (3-1)3^2 \\
   &=3^0 + 3^1 - 3^2 + 3^3\\
   &=11\bar{1}1,
\end{aligned}
 \end{equation*}
as desired.

\subsubsection{An application of signed-digit representation}

An interesting example of how the balanced form of signed-digit representation was used in the past can be found in the work of Fibonacci. With his \emph{Liber Abaci}, Fibonacci spread Hindu-Arabic numbers in Europe together with practical applications, generally of a commercial nature. Zero and the positional number system developed in the Old World on the back of trade and bookkeeping, while in America the calendar was the driving force, and we argue that this difference was to prove crucial. 

The problem is presented as follows: 
\begin{quote}
A certain man in his trade had four weights with which he could weigh integral pounds from one up to 40; it is sought how many pounds was each weight. 
\end{quote}
Fibonacci provides a solution, stating that the four weights are 1lb, 3lb, 9lb and 27lb respectively. 
Clearly, each weight can be used on either side of the balance. We can then give a weight three possible values as a digit: -1 if it is used on the pan with the unknown weight on, 1 if used on the other pan and 0 if the weight is not used at all. 
The most natural choice for a counting system is then the balanced ternary system\footnote{It is interesting that this balanced form of representation minimizes the number of carries in addition, at least when the base is an odd prime \cite{alon2013}.}
In this way, the highest number we can express is 40 (written as 1111 = 27 + 9 + 3 + 1) and the lowest is -40 ($\bar{1}\bar{1}\bar{1}\bar{1}$). Every number within the range -40 to 40 has a unique representation. 
Observe that, from a practical point of view, the weights cannot be negative, so the system is useful only for representing one half of the integer set, that is, the positive ones from 0 to 40.

\subsubsection{The shifting property of digits}

There is another way to think about the foregoing scales problem: 4 weights with 3 possible positions each give rise to $3^4 = 81$ combinations. On the scale, these combinations would read as the 81 numbers from -40 to 40. 
If we were to use a standard ternary system with digits 0, 1, 2,  we could still express exactly 81 consecutive numbers; however these numbers would go from 0 to $2222_{3} = 80_{10}$ (recall the subscript refers to the base).
This idea can easily be generalized: if all digits are shifted by the same quantity, the representation remains non-redundant but the interval of represented numbers is also shifted.
Note that our reasoning proves not only the uniqueness of the representation, but its converse as well; every number is guaranteed to be representable. When using strings of a given length, if we shift all the digits, the interval of representable numbers moves up or down accordingly, but leaves no gaps.\footnote{Note, however that, as in the case of the Fibonacci weights (see below), the useful represented numbers are the positive set. In power systems, negative numbers are externally defined by a minus sign. We can interpret the external minus sign as the additional possibility of making all signs of all digits negative for a given represented number.}

Consider an ordinary base-$b$ system and a string of digits of length $k$, 
$$d_{k-1}d_{k-2}\ldots d_0 = d_{k-1}b^{k-1} + d_{k-2}b^{k-2} \ldots + d_0 b^0,$$  $0 \leq d_i \leq b-1 \quad \forall i.$
Clearly, the smallest representable number is 0 and the largest is 
$$(b-1)(b^{k-1} + \ldots  + b^0) =( b-1) \left(\frac{b^k-1}{b-1} \right) = b^k -1,$$
thus defining an interval of $b^k$ numbers. Again, keeping in mind that we have $b$ possibilities for $k$ positions, the number of combinations is indeed $b^k$.
Suppose now that all the digits are shifted by a quantity $s \in \mathbb{Z}$, i.e., 
$$s \leq d_i \leq b-1+s \quad \forall i.$$
As a result, the lower bound of the interval of representable numbers is 
$$s(b^{k-1} + \ldots  + b^0) = s\left(\frac{b^k-1}{b-1} \right)$$
and similarly the upper bound is 
$$s(b^{k-1} + \ldots + b^0) = (b-s-1)\left(\frac{b^k-1}{b-1} \right).$$
Thus, a shift by $s$ in the digits range results in a shift by $s (b^k-1)/(b-1)$ in the range of represented numbers.

\subsection{Non-power Number Representation Systems}

As we have noted, there are two ways to obtain redundant representation systems. Having seen the signed-digit case, we now move to the other case, that of a non-power representation system. This means that  instead of having a given base or radix,
the positional weights are numbers of a sequence that grows more slowly than the powers of some number. If this is the case, then every number can be represented and generally has more than one expression within the system.\footnote{With a sequence of positional weights growing more rapidly than a power and a fixed set of digits, the positional representation system has gaps, that is, some integer numbers cannot be represented. This is overcome in the factorial number system \cite{knuth}.
}

In the case of a binary system,
we can pick any sequence of numbers that grows more slowly than do powers of two. A famous choice 
 is to use the Fibonacci numbers \cite{fraenkel,butler}.
Fibonacci numbers are the elements of a sequence where each number $F_n$ is the sum of the previous two; $F_n = F_{n-1} + F_{n-2}$, with $F_1 = F_2 = 1$.
The first terms of the sequence are then 1, 1, 2, 3, 5, 8, 13, 21, $\ldots$. These grow more slowly than powers of two.
If one uses these as positional weights, and 0, 1 as digits, it can be seen that every number is representable. Moreover, this representation is unique provided that there are no two consecutive 1's \cite{zeck}.
In fact, since every term is the sum of the previous two, it is easy to see that any string of the form $\ldots011\ldots$ can be replaced with $\ldots100\ldots$ .

Another, more ancient example is a Maya calendrical counting system, the Long Count, which is the subject of this work.

\bibliographystyle{unsrt}   
\bibliography{non-power}

\end{document}